\documentclass[12pt]{article}
 \usepackage{amssymb}
\newtheorem{theorem}{Theorem}

\newtheorem{proposition}[theorem]{Proposition}
\newtheorem{lemma}[theorem]{Lemma}
\newtheorem{corollary}[theorem]{Corollary}
\newtheorem{remark}[theorem]{Remark}
\newtheorem{remarks}[theorem]{Remarks}

\newcommand{\R}{\mathbb{R}}
\newcommand{\Ee}{\mathbb {E}}
\newcommand{\Q}{\mathbb{Q}}

\newcommand{\Sf}{\mathbb{S}}

\newcommand{\Le}{\mathbb{L}}

\newcommand{\Hy}{\mathbb{H}}

\newcommand{\po}{{\hspace*{-1ex}}{\bf .  }}

\def\<{\langle}

\def\>{\rangle}

\def\va{\varphi}

\def\e{\epsilon}
\def\d{\partial}
\def\bea{\begin{eqnarray*} }
\def\eea{\end{eqnarray*} }
\def\be{\begin{equation} }
\def\ee{\end{equation} }

\def\proof{\noindent{\it Proof: }}
\def\qed{\ifhmode\unskip\nobreak\fi\ifmmode\ifinner
\else\hskip5 pt \fi\fi\hbox{\hskip5 pt \vrule width4 pt
height6 pt  depth1.5 pt
\hskip 1pt }}
\begin{document}

\title{ On a class of hypersurfaces 
in  $\Sf^n\times \R$ and $\Hy^n\times \R$.}

\author { Ruy Tojeiro}
\date{}
\maketitle
\addtocounter{equation}{1}

\begin{abstract} We give a complete description of all hypersurfaces of  the product spaces $\Sf^n\times \R$ and $\Hy^n\times \R$ that have flat normal bundle when regarded as submanifolds with codimension two of the underlying flat spaces $\R^{n+2}\supset \Sf^n\times \R$ and $\Le^{n+2}\supset \Hy^n\times \R$. 
We prove that any such hypersurface in $\Sf^n\times \R$ (respectively, $\Hy^n\times \R$) can be  constructed by means of  a family of parallel hypersurfaces in $\Sf^n$ (respectively, $\Hy^n$) and a smooth function of one variable. Then we show that constant mean curvature hypersurfaces in this class are given in terms of an isoparametric family in the base space and a solution of a certain   ODE. For minimal hypersurfaces such solution is explicitly determined in terms of the mean curvature function of the isoparametric family. As another consequence  of our general result,   we  classify the constant angle hypersurfaces of $\Sf^n\times \R$ and $\Hy^n\times \R$, that is, hypersurfaces  with the property that its unit normal vector field makes a constant angle with the unit vector field spanning the second factor $\R$. This extends previous results by Dillen,  Fastenakels,  Van der Veken, Vrancken and Munteanu for surfaces in $\Sf^2\times \R$ and $\Hy^2\times \R$. 
 Our method  also yields   a classification of all Euclidean hypersurfaces with the property that the tangent component of a constant vector field in the ambient space is a principal direction, in particular of all  Euclidean hypersurfaces whose unit normal vector field makes a constant angle with a fixed direction.
 \end{abstract}
\medskip
\section[Introduction]{Introduction} The study of hypersurfaces of  the product spaces $\Sf^n\times \R$ and $\Hy^n\times \R$ has attracted the attention of several geometers  in the last years. Here $\Sf^n$ and $\Hy^n$ denote the sphere and hyperbolic space of dimension $n$, respectively. A natural class of such hypersurfaces consists of those which have flat normal bundle when regarded as submanifolds with codimension two of the underlying flat spaces $\R^{n+2}\supset \Sf^n\times \R$ and $\Le^{n+2}\supset \Hy^n\times \R$, where $\R^{n+2}$ and $\Le^{n+2}$ denote the Euclidean and Lorentzian spaces of dimension $(n+2)$, respectively. Surfaces in $\Sf^2\times \R$ with this property  have been recently classified in \cite{dillen}. The class includes, in particular, all rotational hypersurfaces \cite{veken} and all hypersurfaces with constant sectional curvature and dimension $n\geq 3$ \cite{manfio}. It also contains all {\em constant angle hypersurfaces\/}, that is, hypersurfaces  with the property that its unit normal vector field makes a constant angle with the unit vector field spanning the second factor $\R$. Constant angle surfaces in $\Sf^2\times \R$ and $\Hy^2\times \R$ have been completely described in \cite{vrancken} and \cite{munteanu}, respectively. A similar description was given in \cite{nistor}  of surfaces in $\R^3$ whose unit normal vector field makes a constant angle with a fixed direction.

In order to state our results, let 
 $\Q_\e^n$ denote either  $\Sf^n$, $\R^n$ or  $\Hy^n$, according as $\e=1$, $\e=0$ or $\e=-1$, respectively. Given a hypersurface  $f\colon\,M^n\to \Q_\e^n\times\R$, 
 let $N$ be a unit normal
vector field and let  $\frac{\partial}{\partial t}$ be a unit
vector field tangent to the second factor. 
Then, a vector field  $T$  and a smooth function $\nu$ on $M^n$ are defined by
\be\label{eq:ddt} \frac{\partial}{\partial t}=f_* T+\nu N.\ee

Our first theorem classifies hypersurfaces $f\colon\,M^n\to \Q_\e^n\times\R$ for which $T$ is everywhere a principal direction. Trivial examples are products $M^{n-1}\times \R$, where $M^{n-1}$ is a hypersurface of $\Q_\e^n$, which correspond to the case in which the {\em angle function\/} $\nu$ vanishes identically. More interesting examples arise as follows. Let   $g\colon\, M^{n-1}\to \Q_\e^n$  be a hypersurface and let $g_s \colon\, M^{n-1}\to \Q_\e^n$ be the family of its parallel hypersurfaces, that is, 
\be\label{eq:parallel}g_s(x)=C_\e(s)g(x)+S_\e(s)N(x),\ee
where $N$ is a unit normal vector field to $g$, 
$$
C_\e(s)=\left\{\begin{array}{l}
\cos(s), \,\,\,\mbox{if}\,\,\e=1
\vspace{1.5ex}\\
1, \,\,\,\mbox{if}\,\,\e=0
\vspace{1.5ex}\\
\cosh(s), \,\,\,\mbox{if}\,\,\e=-1
\end{array}\right.\,\,\,\,\,\,\,\,\mbox{and}\,\,\,\,\,\,\,\,\,\,\,
S_\e(s)=\left\{\begin{array}{l}
\sin(s), \,\,\,\mbox{if}\,\,\e=1
\vspace{1.5ex}\\
s, \,\,\,\mbox{if}\,\,\e=0
\vspace{1.5ex}\\
\sinh(s), \,\,\,\mbox{if}\,\,\e=-1.
\end{array}\right.
$$

 Define 
$$f\colon\, M^n:=M^{n-1}\times \R\to \Q_\e^n\times \R$$  by
\be\label{eq:constantangle}f(x,s)=g_s(x)+a(s)\frac{\partial}{\partial t}\ee
for some smooth function $a\colon\, I\to \R$ with nowhere vanishing derivative. 

\begin{theorem}\po\label{thm:main} The map $f$ defines, at  regular points,   a hypersurface that has $T$ as a principal direction.
Conversely,  any hypersurface $f\colon\, M^n\to \Q_\e^n\times \R$  with nowhere vanishing angle function that has $T$ as a principal direction is locally given in this way.
\end{theorem}

Besides  (open subsets of) slices $\Q_\e^n\times \{t\}$ and products $(M^{n-1}\subset \Q_\e^n)\times \R$, the hypersurfaces   in Theorem \ref{thm:main} turn out to be precisely  the ones that have flat normal bundle when regarded as submanifolds of  $\Ee^{n+2}$ (see Proposition \ref{lem:Rperpzero} in Section $2$). Here   $\e\in \{-1, 1\}$ and $\Ee^{n+2}$ stands for either Euclidean space $\R^{n+2}$ or Lorentzian space $\Le^{n+2}$, according as $\e=1$ or $\e=-1$, respectively.

As a consequence of Theorem \ref{thm:main}, we get the following complete description of all hypersurfaces $f:M^n\to\Q_\e^n\times\R$ with constant angle function, called {\em constant angle hypersurfaces.\/}

\begin{corollary}\po\label{cor:constantangle} Let $f$ be given by (\ref{eq:constantangle}) with $a(s)=As$ for some $A\neq 0$. Then $f$ is a constant angle hypersurface.  Conversely,  any constant angle hypersurface $f\colon\, M^n\to \Q_\e^n\times \R$ is  either an open subset of a  slice $\Q_\e^n\times \{t_0\}$ for some $t_0\in \R$, an open subset of a product $M^{n-1}\times \R$, where $M^{n-1}$ is a hypersurface of $\Q_\e^n$, or it is locally given  in this way.
\end{corollary}

Our next result characterizes constant mean curvature  hypersurfaces of $\Q_\e^n\times \R$ within the class  of those which have  $T$ as a principal direction.

\begin{theorem}\po\label{thm:cmc} Let $g\colon\, M^{n-1}\to \Q_\e^n$ be an isoparametric hypersurface, and let $H(s)$ be the (constant) mean curvature of its parallel hypersurface $g_s$. Given $H\in \R$, let $a\colon\,I\to \R$ be a solution of 
\be\label{eq:edo}a''(s)-a'(s)(1+(a'(s))^2)H(s)-H(1+(a'(s))^2)^{3/2}=0\ee
on an open interval $I\subset \R$ such that $g_s$ is an immersion on $M^{n-1}$ for every $s\in I$.  Then $f$ is a hypersurface with  constant mean curvature $H$ that has $T$ as a principal direction. Conversely, any hypersurface $f\colon\, M^{n}\to \Q_\e^n\times \R$  with  nowhere vanishing angle function and constant mean curvature $H$ that has  $T$ as a principal direction  is 
 locally given in this way.
\end{theorem}

In the minimal case, the ODE (\ref{eq:edo}) can be explicitly solved in terms of the mean curvature function $H(s)$ of the isoparametric family: 

\begin{corollary}\po\label{cor:min} Let $g\colon\, M^{n-1}\to \Q_\e^n$ be an isoparametric hypersurface, and let $H(s)$ be the (constant) mean curvature of its parallel hypersurface $g_s$. Let $I\subset \R$ be an open interval such that  $g_s$ is an immersion on $M^{n-1}$ for every $s\in I$. Given $a_0, h_0\in \R$ with $0<h_0<1$, define  $a\colon\,I\to \R$ by
\be\label{eq:a}a(s)=a_0\pm \int_0^s\sqrt{\frac{h(t)}{1-h(t)}}\,\,dt,\,\,\,\,\mbox{with}\,\,\,\,h(t)=h_0\exp  \left(2\int_0^t H(\tau)\,d\tau\right).\ee
Then  $f\colon\, M^{n-1}\times I\to \Q_\e^n\times \R$  given by (\ref{eq:constantangle}) is a minimal hypersurface that has  $T$ as a principal direction. Conversely, any minimal hypersurface $f\colon\, M^{n}\to \Q_\e^n\times \R$  with  nowhere vanishing angle function that has  $T$ as a principal direction is 
locally given in this way.
\end{corollary}

\section[Preliminaries]{Preliminaries}

 Given a hypersurface  $f\colon\,M^n\to \Q_\e^n\times\R\subset \Ee^{n+2}$ with a unit  normal
vector field $N$, let $A$ be the shape operator of
$f$ with respect to $N$ and let $\nabla$ be the Levi-Civita connection  of $M^n$. Using that
$\frac{\partial}{\partial t}$ is parallel in $\Q_\e^n\times\R$, we obtain by differentiating (\ref{eq:ddt}) that 
\begin{eqnarray}\label{eq:NablaT}
\nabla_XT=\nu AX
\end{eqnarray}
and 
\begin{eqnarray}\label{eq:DerivadaNu}
X(\nu)=-\langle AX,T\rangle,
\end{eqnarray}
for all $X\in TM$. 

Another fact that we will need in the proof of Theorem \ref{thm:main}  is that the vector field $T$ is a gradient vector field. Namely, it is  the gradient of the height function $h=\<f,\frac{\partial}{\partial t}\>$.

As a final observation in this short section, let $\xi$ denote the outward pointing  unit normal vector field to $\Q_\e^n\times\R$ along $f$ and  let $A_\xi$ be  the corresponding shape operator. Then, it is easily seen that 
\be\label{eq:eigenvectorAxi}
A_\xi T=-\nu^2T\,\,\,\,\mbox{and}\,\,\,\,A_\xi X=-X\,\,\,\mbox{for}\,\,X\in \{T\}^\perp.
\ee
This leads to the following characterization of flatness of the normal bundle of  $f$ when regarded as an isometric immersion into $\Ee^{n+2}$, first proved in \cite{dillen} for surfaces in $\Sf^2\times \R$.

\begin{proposition}\label{lem:Rperpzero} Let $f:M^n\to\Q_\e^n\times\R\subset \Ee^{n+2}$, $\e\in \{-1, 1\}$, be a hypersurface. Suppose that $T$ does not vanish at $x\in M^n$. Then $f$ has flat normal bundle at $x$  as an isometric immersion into $\Ee^{n+2}$  if and only if $T$ is a principal direction of $f$ at $x$.
\end{proposition}
\proof
By the Ricci equation, $f$ has flat normal bundle (as an isometric immersion into $\Ee^{n+2}$) if and only if $A$ commutes with $A_\xi$. This is the case if and only if the eigenspaces of $A_\xi$ are invariant by $A$, which by (\ref{eq:eigenvectorAxi})  is equivalent to $T$ being an eigenvector of $A$.\vspace{2ex}\qed

\section[The proofs]{The proofs}

\noindent {\em Proof of Theorem \ref{thm:main}:}
We have 
$$f_*X={g_s}_*X,\,\,\,\,\,\mbox{for any}\,\, X\in TM^{n-1},$$
and 
$$f_*\frac{\d}{\d s}=N_s+a'(s)\frac{\d}{\d t},$$
where \be\label{eq:normalparallel}N_s(x)=-\e S_\e(s)g(x)+C_\e(s)N(x).\ee
Therefore, a point $(x,s)\in M^{n-1}\times \R$ is regular for $f$ if and only if $g_s$ is regular at $x$, in which case
$N_s(x)$ is a unit normal vector  to $g_s$ at $x$ and 
\be\label{eq:eta}\eta(x,s)=-\frac{a'(s)}{b(s)}N_s(x)+\frac{1}{b(s)}\frac{\d}{\d t},\,\,\,\,\,\,\mbox{with}\,\,\,\,b(s)=\sqrt{1+a'(s)^2},\ee
is a unit normal vector  to $f$ at $(x,s)$.
Notice that $\<f_*X,f_*\frac{\d}{\d s}\>=0$ for any $X\in TM^{n-1}$.
We have
$$\tilde{\nabla}_{\d/\d s}\eta=-\left(\frac{a'(s)}{b(s)}\right)'N_s+\left(\frac{1}{b(s)}\right)'\frac{\d}{\d t}+\e\frac{a'(s)}{b(s)}g_s,$$
where $\tilde{\nabla}$ stands for the derivative in $\Ee^{n+2}$. 
Then
$$\< \tilde{\nabla}_{\d/\d s}\eta,f_*X\>=\< \tilde{\nabla}_{\d/\d s}\eta,{g_s}_*X\>=0,\,\,\,\,\mbox{for any $X\in TM^{n-1}$},$$
which shows that ${\d/\d s}$ is a principal direction of $f$. Moreover, using  that 
\be\label{eq:anglefunction}\nu=\<\eta, \frac{\d}{\d t}\>=\frac{1}{b(s)},\ee
we obtain
$$f_*T=\frac{\d}{\d t}-\nu\eta=\frac{a'(s)}{b^2(s)}f_*\frac{\d}{\d s},$$
hence
$$T=\frac{a'(s)}{b^2(s)}\frac{\d}{\d s}.$$
Therefore  $T$ is  a principal direction of $f$.

We now prove the converse. Since $T$ is a gradient vector field,
 the orthogonal distribution $\{T\}^\perp$ is integrable. Hence, there exists locally a diffeomorphism $\psi\colon\,M^{n-1}\times I\to M^n$, where $I$ is an open interval containing $0$, such that $\psi(x,\cdot)\colon\, I\to M^n$ are  integral curves of $T$ for any $x\in M^{n-1}$ and $\psi(\cdot, s)\colon\, M^{n-1}\to M^n$ are  leaves of $\{T\}^\perp$ for any $s\in I$. In particular, 
 $\psi_*X\in \{T\}^\perp$ for any $X\in TM^{n-1}$. 
Set $F= f\circ \psi$.  Then 
$$X\<F, \frac{\d}{\d t}\>=\<f_* \psi_*X, \frac{\d}{\d t}\>=\<\psi_*X, T\>=0$$
for any $X\in TM^{n-1}$. Thus $\<F(x,s), \frac{\d}{\d t}\>=\rho(s)$ for some smooth function $\rho$  on $I$.

We claim that $\Pi_1\circ F(x,\cdot)\colon\, I\to \Q_\e^n$   is a  pre-geodesic of $\Q_\e^n$ for any $x\in M^{n-1}$, where  $\Pi_1\colon\, \Q_\e^n \times \R\to \Q_\e^n$ is the canonical projection, that is, the arclength reparametrization of $\Pi_1\circ F(x,\cdot)$ is a geodesic of $\Q_\e^n$.  In other words, $\alpha:=\Pi_1\circ f\circ \gamma$ is a pre-geodesic of $\Q_\e^n$ for any  integral curve $\gamma$ of $T$. 

First notice that, since   $T$ is a principal direction of $f$, it follows from (\ref{eq:DerivadaNu})  that $X(\nu)=0$ for any $X\in \{T\}^\perp$, hence also  $X(\|T\|)=0$ for any $X\in \{T\}^\perp$,  for $\|T\|^2+\nu^2=1$. Then, the  following general fact implies that
$\gamma$ is a pre-geodesic of $M^n$.

\begin{lemma}\po \label{le:pregeodesic} Let $T$ be a gradient vector field on a Riemannian manifold $M^n$. Assume that $\|T\|$ is constant along $\{T\}^\perp$. Then the integral curves of $T$ are pre-geodesics of $M^n$.
\end{lemma}
\proof Since $T$ is a gradient vector field, we have
$$\<\nabla_XT,Y\>=\<\nabla_YT,X\>$$
for all $X,Y\in TM$. Therefore, for any $X\in \{T\}^\perp$ we obtain that
$$\<\nabla_T T,X\>=\<\nabla_X T,T\>=\frac{1}{2}X(\|T\|^2)=0.\qed$$

Now observe that the velocity vector of $\alpha$  is
$f_*T-\<f_*T,{\d}/{\d t}\>{\d}/{\d t}$,
whose length is $\lambda=\|T\|\nu$. Therefore, all we need  to prove is that 
$$\tilde{\nabla}_T (f_*(\lambda^{-1}T)-\<f_*(\lambda^{-1}T),{\d}/{\d t}\>{\d}/{\d t})$$
lies in the direction of the normal vector field $\xi$ to $\Q_\e^n\times \R$ along $f$. We have
\be\label{eq:1} \tilde{\nabla}_T f_*(\lambda^{-1}T)=T({\nu}^{-1})f_*(\hat{T})+{\nu}^{-1}\tilde{\nabla}_T f_*\hat T,\ee
where $\hat T=T/\|T\|$.
Now, since ${\nabla}_T \hat{T}=0$ by Lemma \ref{le:pregeodesic}, we have using (\ref{eq:NablaT}) and (\ref{eq:eigenvectorAxi}) that
\be\label{eq:2}\tilde{\nabla}_T f_*\hat{T}=\<AT,\hat{T}\>N+\<A_\xi T,\hat{T}\>\xi=\nu^{-1}T(\|T\|)N-\nu^2\|T\|\xi.\ee
Then, from (\ref{eq:1}),  (\ref{eq:2}) and 
$$T(\nu^{-1})=-\nu^{-2}{T(\nu)}=-(1/2)\nu^{-3}{T(\nu^2)}={(1/2)\nu^{-3}}T(\|T\|^2)={\nu^{-3}}{\|T\|T(\|T\|)},$$
we obtain 
\be\label{eq:3}\tilde{\nabla}_T f_*(\lambda^{-1}T)=\nu^{-3}\|T\|T(\|T\|)f_*\hat{T}+\nu^{-2}T(\|T\|)N-\nu\|T\|\xi=\nu^{-3}T(\|T\|)\frac{\d}{\d t}-\nu\|T\|\xi.\ee
On the other hand, we have $\<f_*(\lambda^{-1}T),\frac{\d}{\d t}\>=\nu^{-1}\|T\|,$
and 
\be\label{eq:4} T(\nu^{-1}\|T\|)=T(\nu^{-1})\|T\|+\nu^{-1}T(\|T\|)=\nu^{-3}T(\|T\|)(\|T\|^2+\nu^2)=\nu^{-3}T(\|T\|).\ee
It follows from (\ref{eq:3}) and (\ref{eq:4})
that
$$\tilde{\nabla}_T (f_*(\lambda^{-1}T)-\<f_*(\lambda^{-1}T),{\d}/{\d t}\>{\d}/{\d t})=-\nu\|T\|\xi,$$
which proves the claim.\vspace{1ex}

Now, since $\|T\|$ and $\nu$ are constant along $\{T\}^\perp$,   there exists a smooth function $r\colon\,I\to \R$ such that $\|T\|\nu\circ \psi(x,s)=r(s)$ for all $(x,s)\in M^{n-1}\times I$.  Define $g\colon\,M^{n-1}\to \Q_\e^n$ by $g=\Pi_1\circ F$, and let $g_s$ be the family of parallel hypersurfaces to $g$. Set $\va(s)=\int_{s_0}^s r(\sigma)d\sigma$, $a=\rho\circ \va^{-1}$ and $\tilde \psi(x,s)=\psi(x, \va^{-1}(s))$ for $(x,s)\in M^{n-1}\times J$, with $s_0\in I$ and $J=\va(I)$. By the claim, we have that  
 $$ f\circ \tilde \psi(x,s)=g_s(x)+a(s)\frac{\d}{\d t}\,\,\,\mbox{for any $(x,s)\in M^{n-1}\times J$}. \,\,\,\qed$$

\begin{remarks}{\em $(i)$ We have seen that a point $(x,s)\in M^{n-1}\times \R$ is regular for $f$ if and only if $g_s$ is regular at $x$. Let us discuss when the latter occurs. Let $\lambda_1,\ldots, \lambda_m$ be the distinct principal curvatures of $g$, excluding $0$ if $\e=0$ and those with absolute value less than or equal to $1$ if $\e=-1$. 
For $1\leq i\leq m$, write
\be\label{lambdas}
\lambda_i=\left\{\begin{array}{l}\cot \theta_i,\,\,\,\,0<\theta_i<\pi,\,\,\,\mbox{if}\,\,\,\e=1,\vspace{1.5ex}\\
\coth \theta_i,\,\,\,\,\theta_i\neq 0,\,\,\,\mbox{if}\,\,\,\e=-1,\vspace{1.5ex}\\
1/\theta_i,\,\,\,\,\theta_i\neq 0,\,\,\,\mbox{if}\,\,\,\e=0,\vspace{1.5ex}\\
\end{array}\right.\ee 
where the $\theta_i$ form an increasing sequence. If $X$ is in the  eigenspace of the shape operator $A_N$ corresponding to the principal curvature $\lambda_i$, $1\leq i \leq m$, we have
\be\label{gs*}
{g_s}_*X=\left\{\begin{array}{l}{\displaystyle \frac{\sin(\theta_i-s)}{\sin \theta_i}X,\,\,\,\mbox{if}\,\,\,\e=1,}\vspace{1.5ex}\\
{\displaystyle \frac{\sinh(\theta_i-s)}{\sinh \theta_i}X,\,\,\,\mbox{if}\,\,\,\e=-1,}\vspace{1.5ex}\\
{\displaystyle \frac{\theta_i-s}{\theta_i}X,\,\,\,\mbox{if}\,\,\,\e=0.}
\end{array}\right.\ee
 Thus, for $\e=0$ and $\e=-1$ (respectively, $\e=1$), $g_s$ is an immersion at $x$ if and only if $s\neq \theta_i(x)$ (respectively, $s\neq \theta_i(x)\mbox{(mod $\pi$)}$) for any $1\leq i \leq m$. If $\e=0$ (respectively,  $\e=-1$), let $\theta_+$  be the least  of the $\theta_i$ that is  greater than $0$ (respectively, $1$), and let  $\theta_-$  be the greater  of the $\theta_i$ that is  less than $0$ (respectively, $-1$). Set
\be\label{U}
U:=\left\{\begin{array}{l}{\displaystyle \{(x,s)\in M^{n-1}\times \R\,:\,s\in (\theta_m(x)-\pi, \theta_1(x))\},\,\,\,\mbox{if}\,\,\,\e=1,}\vspace{1.5ex}\\
{\displaystyle \{(x,s)\in M^{n-1}\times \R\,:\,s\in (\theta_-(x), \theta_+(x))\},\,\,\,\mbox{if}\,\,\,\e=0 \,\,\,\mbox{or}\,\,\,\e=-1.}
\end{array}\right.\ee
In any case, if  $V\subset M^{n-1}$ is an open subset and $I$ is an open interval containing $0$ such that $V\times I\subset U$, then $g_s$ is an immersion on $V$ for every $s\in I$, and hence $f$ is an immersion on $V\times I$. In particular, if $g$ is an isoparametric hypersurface, one can take $V=M^{n-1}$ and $I= (\theta_m-\pi, \theta_1)$ if $\e=1$ and $I=(\theta_-, \theta_+)$ if $\e=0$ or $\e=-1$.\vspace{2ex}\\
$(ii)$ The hypersurface $f$ given by (\ref{eq:constantangle}) has a nice geometric description in terms of $g\colon\,M^{n-1}\to \Q_\e^{n}$. Assume first that $\e=\pm 1$. Regarding $g$ as an isometric immersion into $\Ee^{n+2}$, its normal space at each point $x\in M^{n-1}$ is a Lorentzian or Riemannian vector space of dimension $3$, according as $\e=-1$ or $\e=1$, respectively, which is spanned by the position vector $g(x)$, the normal vector  $N(x)$ to $g$ in $\Q_\e^{n}$ at $x$ and the constant vector $\d/\d t$.  Notice that these give rise to parallel vector fields  in the normal connection of $g$. For a fixed $x\in M^{n-1}$, we can regard $f(x,s)=C_\e(s)g(x)+S_\e(s)N(x)+a(s)\d/\d t$ as a curve in a cylinder $\Q_\e^1\times \R$ with axis $\d/\d t$ contained in the normal space of $g$ at $x$. Thus, the immersion $f$ is generated by parallel transporting such curve in the normal connection of $g$. Moreover, constant angle hypersurfaces correspond to the case in which such curve is a helix in $\Q_\e^1\times \R$. A similar description holds for $\e=0$. Hypersurfaces of $\R^{n+1}=\R^n\times \R$ with the property that the tangent component $T$ of the constant vector field $\d/\d t$ spanning the second factor is a principal direction are generated by parallel transporting a curve in the normal space in $\R^{n+1}$ of a hypersurface $g$ of $\R^n$. Constant angle hypersurfaces arise in the particular case in which such a curve is a straight line.}
\end{remarks}

\noindent {\em Proof of Corollary \ref{cor:constantangle}:$\,$} The direct statement follows from (\ref{eq:anglefunction}). If $\nu=1$ or $\nu=0$, it is easily seen that   $f(M^n)$ is  an open subset of a  slice $\Q_\e^n\times \{t_0\}$,  $t_0\in \R$, or  an open subset of a product $M^{n-1}\times \R$, where $M^{n-1}$ is a hypersurface of $\Q_\e^n$, respectively. Otherwise,  it follows from (\ref{eq:DerivadaNu}) that the vector field $T$ is a principal direction. By Theorem \ref{thm:main},   the hypersurface $f$ is  locally given by (\ref{eq:constantangle}). Finally, since $\nu$  is given by 
(\ref{eq:anglefunction}), the fact that it is constant 
 implies that the function $a(s)$ in  (\ref{eq:constantangle}) is linear (with nowhere vanishing derivative), hence we may assume that $a(s)=As$ for some $A\neq 0$.
\vspace{2ex}\qed\\
\noindent {\em Proof of Theorem \ref{thm:cmc}:} 
Let $f\colon\, M^{n-1}\times \R\to \Q_\e^n\times \R$ be given by (\ref{eq:constantangle}). We have
$$\tilde{\nabla}_X\eta= \frac{a'(s)}{b(s)}{g_s}_*A^s X\,\,\,\,\,\mbox{for any}\,\, X\in TM^{n-1},$$
hence the shape operator $A_\eta$ satisfies
$$A_\eta X=-\frac{a'(s)}{b(s)}A^s X,\,\,\,\,\,\mbox{for any}\,\, X\in TM^{n-1}.$$
On the other hand,
$$\< \tilde{\nabla}_{\d/\d s}\eta,f_*{\d/\d s}\>=-\left(\frac{a'(s)}{b(s)}\right)'+a'(s)\left(\frac{1}{b(s)}\right)'=-\frac{a''(s)}{b(s)},$$
hence the principal curvature in the ${\d/\d s}$- direction is $a''(s)/b^3(s)$. It follows that the (non normalized)
mean curvature function of $f$ is given by
\be\label{eq:edo2}H=-\frac{a'(s)}{b(s)}H_s+\frac{a''(s)}{b^3(s)}=\frac{-a'(s)(1+(a'(s))^2)H_s+a''(s)}{(1+(a'(s))^2)^{3/2}},\ee
where $H_s$ denotes the mean curvature function of $g_s$. The conclusion now follows from the fact that $H_s$ is  constant on $V$ (that is, it depends only on $s$) if and only if $g$ is an isoparametric hypersurface of $\Q_\e^n$ (see  \cite{cecil}, Theorem $5.8$, p. $272$). \vspace{2ex}\qed\\
\noindent {\em Proof of Corollary \ref{cor:min}:$\,$} Let us prove the converse.  By Theorem \ref{thm:cmc}, if  $f\colon\, M^{n}\to \Q_\e^n\times \R$  is a minimal hypersurface with  nowhere vanishing angle function that has  $T$ as a principal direction, then it is  
locally given by (\ref{eq:constantangle}) for some solution $a\colon\,I\to \R$ of (\ref{eq:edo}) with $H=0$, that is, 
$$a''(s)-a'(s)(1+(a'(s))^2)H(s)=0.$$
Without loss of generality, we may assume that  $a'(s)>0$ for all $s\in I$.  Then
$$\log\frac{a'(s)}{\sqrt{1+(a'(s))^2}}=\log\frac{a'(0)}{\sqrt{1+(a'(0))^2}}+\int_0^s H(\tau)\,d\tau,$$
 hence
$$\frac{(a'(s))^2}{{1+(a'(s))^2}}=h(s):=\frac{(a'(0))^2}{1+(a'(0))^2}\exp  \left(2\int_0^s H(\tau)\,d\tau\right).$$
It follows that $0<h(s)<1$ for all $s\in I$ and 
$$a(s)=a(0)+\int_0^s\sqrt{\frac{h(t)}{1-h(t)}}\,\,dt.$$
Had we assumed that $a'(s)<0$ for all $s\in I$, we would have obtained the same expression for $a$, but with the minus sign in (\ref{eq:a}). The conclusion follows by taking $a_0=a(0)$ and $h_0=h(0)={(a'(0))^2}/({1+(a'(0))^2})$. The direct statement is now clear.\qed

\begin{remark}{\em One can check that for $\e=0$ the minimal hypersurfaces given by Corollary~\ref{cor:min} are the minimal $n$-dimensional catenoids described in \cite{dcd}.}
\end{remark}

\end{document}